\documentclass[12pt]{amsart}
\usepackage{amssymb,amsmath}
\numberwithin{equation}{section}
\def\p{\partial}

\def\1#1{\overline{#1}}
\def\2#1{\widetilde{#1}}
\def\3#1{\widehat{#1}}
\def\4#1{\mathbb{#1}}
\def\5#1{\frak{#1}}
\def\6#1{{\mathcal{#1}}}
\def\C{{\4C}}

\def\B{\Bbb B}

\def\phi{\varphi}

\def\D{\Delta}

\newtheorem{Thm}{Theorem}[section]
\newtheorem{Cor}[Thm]{Corollary}
\newtheorem{Pro}[Thm]{Proposition}
\newtheorem{Lem}[Thm]{Lemma}
\theoremstyle{definition}\newtheorem{Def}[Thm]{Definition}
\theoremstyle{remark}
\newtheorem{Rem}[Thm]{Remark}
\newtheorem{Exa}[Thm]{Example}

\def\Label#1{\label{#1}}
\def\bl{\begin{Lem}}
\def\el{\end{Lem}}
\def\bp{\begin{Pro}}
\def\ep{\end{Pro}}
\def\bt{\begin{Thm}}
\def\et{\end{Thm}}
\def\bc{\begin{Cor}}
\def\ec{\end{Cor}}
\def\bd{\begin{Def}}
\def\ed{\end{Def}}
\def\br{\begin{Rem}}
\def\er{\end{Rem}}
\def\be{\begin{Exa}}
\def\ee{\end{Exa}}
\def\bpf{\begin{proof}}
\def\epf{\end{proof}}
\def\ben{\begin{enumerate}}
\def\een{\end{enumerate}}

\def\1alpha{[\frac1\alpha]}

\def\C{{\Bbb C}}

\numberwithin{equation}{section}

\newtheorem{theorem}{Theorem  }[section]
\newtheorem{definition}[theorem]{Definition }
\newtheorem{lemma}[theorem]{Lemma  }
\newtheorem{proposition}[theorem]{Proposition  }
\newtheorem{corollary}[theorem]{Corollary }
\newtheorem{example}[theorem]{\it Example }

\begin{document}
\title[Holomorphic extension from the sphere to the ball]{Holomorphic extension from the sphere to the ball}
\author[L.~Baracco]
{Luca Baracco}
\address{Dipartimento di Matematica, Universit\`a di Padova, via 
Trieste 63, 35121 Padova, Italy}
\email{baracco@math.unipd.it}
\maketitle
\begin{abstract}
Real analytic functions on the boundary of the sphere which have separate holomorphic extension along the complex lines through a boundary point have holomorphic extension to the ball. This was proved in \cite{B09} by an argument of CR geometry. We give here an elementary proof based on the expansion in holomorphic and antiholomorphic powers.
\end{abstract}
\section{Main result - Statement and proof}
In my note \cite{B09} I stated a general principle of holomorphic extension 
from a convex boundary of $\C^n$ for functions which have separate holomorphic extension along generic $(2n-2)$-parameter families of discs.  In particular, the discs which pass through a fixed boundary point. The conclusion seems to be natural in the framework of CR geometry but less in the classical literature of harmonic analysis. Also, it is in contrast with the case in which the ``center" of the system of discs is an interior point of the ball: a single point does not suffice.  
For this reason, I want to supply here a direct proof which only uses Taylor expansions: moment condition forces coefficients of the antiholomorphic powers to vanish. The technique is specific of the sphere. General convex domains are out of reach unless one brings into play, as I do in \cite{B09}, the extraordinary strength of the Riemann-Lempert Theorem (cf. \cite{L81} and \cite{CL88}). This makes equivalent the family of discs by a point of a convex domain to the lines by a point of the ball.

Let $\B^n$ be the ball in $\C^n$, $\partial \B^n$ the sphere, $z_o$ a point in $\partial\B^n$, and $C^\omega$ the real analytic functions.
\bt
\Label{t1.1}
Let $f$ be a function in $C^\omega(\partial\B^n)$ and suppose that $f$ extends from $\partial\B^n$ along each line  passing through $z_o$. Then $f$ extends holomorphically to $\B^n$.
\et
\bpf
{\bf (a)} We first prove the result for $\B^2$ in $\C^2$. It is not restrictive that $z_o$ is the pole $(0,1)$.
The straight discs through $(0,1)$ can be parametrized  over $a\in\C$ as the sets $D_a$ described by
$$ D_a(\tau)=\left( \frac{\tau -1}{1+|a|^2} a,\frac{\tau -1}{1+|a|^2} +1 \right)\quad \forall\tau\in\bar{\D} .$$
Note that when $|a|>>1$ the disc $D_a$ becomes very close to the complex tangent line to the sphere at the point $z_o$, and moreover $D_a$ lies in a neighborhood of $z_o$.

Since $f\in C^\omega(\p\B^2)$, and $\bar{\p}_{z_2}$ is transverse to $\p\B^2$ at $z_o$,  $f$ can be extended in a neighborhood of $z_o$ holomorphically in $z_2$. We denote  again by $f$ this extension.
We consider the power series development of $f$ at $z_o$
$$ f(z_1,\bar{z}_1,z_2)=\sum^{+\infty}_{l=0} \sum_{h+k+2m=l} b_{h,k,m} z_1^h\bar{z}_1^kz_2^m$$
note that we reordered the terms in a weighted degree (giving  weight $2$  to  $z_2$).
Taking $|a|$ sufficiently big we consider the $N$-momentum on the disc $D_a$:
\begin{align}
G(a,N)&=\int_{\p\D} \tau^N f(D_a(\tau))d\tau 
\\ &=\int_{\p\D}\tau^N \sum^{+\infty}_{l=0} \sum_{h+k+2m=l} b_{h,k,m}
\left(\frac{\tau -1}{1+|a|^2} a\right)^h \left(\overline{\frac{\tau -1}{1+|a|^2} a}\right)^k\left(\frac{\tau -1}{1+|a|^2}\right)^m 
\end{align}
We want to prove that $b_{h,k,m}=0$ whenever $k>0$. To this end, let $l_o$ be the lowest weighted degree such that $b_{h,k,m}\neq 0$ for some $k>0$ and let $k_o$ be the highest degree in $\bar{z}_1$ for which this happens.
We get $G(a,N)=0$ for any $N$ and any $a$, in particular, for $ta$ with $|a|=1$ and $t\to +\infty$. 
Consider the  limit
\begin{align}
 \lim_{t\to +\infty} {G(ta,N)}{t^{l_o}}
&=\lim_{t\to +\infty}\sum_{l=l_o}^{+\infty}\sum_{h+k+2m=l,k>N} 
 2\pi i(-1)^{k+h+N-1}\times
 \\
 &\qquad\times\binom{h+k+m}{k-N-1} a^h\bar{a}^k\left(\frac1{t^2}+|a|^2)^m\right) t^{l_o-l}\frac1{(\frac1{t^2}+|a|^2)^l}\\&=
\sum_{h+k+2m=l_o, k>N}  (2\pi i)(-1)^{h+k+N-1} \binom{h+k+m}{k-N-1} b_{h,k,m} \frac{a^h\bar{a}^k|a|^{2m}}{|a|^{2l_o}}=0,
\end{align}
where we have used the fact  that $\int_{\p\D} \tau^N(\tau-1)^h(\bar{\tau}-1)^k(\tau -1)^m d\tau=\int_{\p\D}(-1)^k\frac{ \tau -1)^{h+k+m}}{\tau^(k-l)} d\tau=(-1)^{h+m+N+1}\binom{h+k+m}{k-1-N}$. 
Now, choosing $N=k_o-1$, we get the following relation on the coefficients $b$'s:
$$ \sum_{h+k_o+2m=l_o}(-1)^{h+m+k_o} \binom{h+k_o+m}{1}b_{h,k_o,m}a^{h+m}\bar{a}^{k_o+m}=0.$$
Writing $a=e^{i\theta}$, we get
$$\sum_{h+k_o+2m=l_o}(-1)^{h+m+k_o} \binom{h+k_o+m}{1}b_{h,k_o,m}e^{i\theta(h-k_o)}=0,$$
which implies $b_{h,k_o,m}=0$ for $h+k_o+2m=l_o$. Therefore, when $k\geq 1$, we have $b_{h,k,m}=0$ for  any weighted degree $l$.
This concludes the proof in dimension $2$.

\noindent
{\bf (b)} We pass to the ball $\B^n$ of general dimension. We still suppose that $z_o$ is the pole $(0,...,1)$. By (a) we know that $f$ extends holomorphically along the slices of $\B^n$ with the 2-dimensional planes through $0$ and $z_o$. 
The different extensions glue to a single, well defined, function on $\B^n$ because, on the complex line $L$ through $0$ and $z_o$ where they all overlap, they have to coincide by the Caucly formula over  $L\cap\partial\B^n$.
Also, by Cauchy-Kowalevski Theorem, $f$ extends to a neighborhood of $z_o$ in $\B^n$ as a  holomorphic  function of $z_n$. 
This must coincide with the extension along the slices but it is, in addition, $C^\omega$.
In any of these slices $f$ is holomorphic in radial directions and therefore, by Forelli's Theorem (cf. for instance \cite{R}), it is holomorphic. Since it is also holomorphic in $z_n$, then it is holomorphic in a neighborhood of $z_o$ in $\B^n$. By real analyticity, the holomorphicity of the extension propagates from a neighborhood of $z_o$ to a neighborhood of the whole $\partial\B^n$. Finally, by Hartogs' Theorem, $f$ is in fact holomorphic in the whole $\B^n$.

\epf

\end{document}